\documentclass[12pt]{article}
\usepackage[latin1]{inputenc} 
\usepackage{graphicx}
\usepackage{amssymb}
\usepackage{fontawesome}
\usepackage{amsthm}
\usepackage{amsmath}
\usepackage{url}
\usepackage{tikz}
\usetikzlibrary{decorations.markings}
\tikzstyle{point}=[circle, draw, inner sep=0pt, minimum size=6pt]
\usepackage{caption}
\theoremstyle{plain}

\theoremstyle{definition} 
\theoremstyle{remark} 

\theoremstyle{definition}

\begin{document}
\begin{center}
\large{\bf{Steiner symmetrization with respect to the Kakutani-Fibonacci sequence of directions}}
\bigskip\bigskip\bigskip

{\large{\bf{Ingrid Carbone and Aljo\v{s}a Vol\v{c}i\v{c}}}}

Universit\`a della Calabria
\bigskip

Dedicated to the memory of our dear friend Paolo Gronchi.
\end{center}
\bigskip

\begin{abstract}
In this paper we will prove that for any planar measurable set of finite measure $M$, its successive Steiner symmetrizations with respect to the Kakutani-Fibonacci sequence of directions converge to the ball $M^*$ centered at the origin and having the same measure.
\end{abstract}

\section {\bf Introduction} 

The paper is motivated by Conjecture 5.2 of [BBGV] which asks whether sequences of directions having low-discrepancy generate a convergent Steiner process. \smallskip

In the attempt to prove the long-standing isoperimetric problem, Jakob Steiner [S] introduced a beautiful geometric construction which is named after him: the Steiner symmetrization. 

Despite of ingenuity of the construction, his solution was incomplete as he did not prove its existence.

Steiner's result was later corrected and extended to higher dimension and more general classes of sets. The first complete proof was obtained by Weierstrass using analytic methods. Due to its significance (and beauty), the isoperimetric problem has been studied over the years. 

The seminal paper of De Giorgi [DG] marked a milestone by proving the isoperimetric property of the ball within the class of measurable sets with finite measure and finite Caccioppoli-De Giorgi perimeter.

A common feature of many papers concerning Steiner symmetrization is the following: starting from a set $E$ (called the seed) of a certain class (convex bodies, compact sets, measurable sets with finite measure), a first symmetrization is taken by choosing, in a clever way, a direction $u_1$ to obtain the set $S_{u_1}E$, called its {\it Steiner symmetrization}, which is ``more round". With a second clever choice of a direction $u_2$ one gets $S_{u_2}S_{u_1}E$, which is even more round. The construction continues selecting directions $u_3\,,u_4, \dots, u_n, \dots$, obtaining a sequence of sets $\{S_{u_n}\dots S_{u_1}E\}$, called {\it Steiner process}. The desired conclusion follows if the process converges (in the appropriate distance) to the ball $E^*$, having the same volume as $E$, centered at the origin $o$. Directions are represented as elements of $S^1_+ =\{u: u=e^{i\vartheta}, \vartheta \in [0,\pi]\}$.

The choice of the $n$-th direction is made step by step, guided by the shape of the $(n-1)$-th set and the decrease (or increase) of an appropriate functional like perimeter, moment of inertia or others. The use of the axiom of choice is inevitable and the proof is not constructive.

As far as the authors know, the first algorithmic construction of a convergent Steiner process, in any dimension, is due to Eggleston [E]. In dimension two it consists in taking two directions $u$ and $v$ making an angle which is an irrational multiple of $\pi$ and alternating the symmetrizations in the two directions. He proved that $(S_uS_v)^nC$ tends to $C^*$ for any convex body $C$.

There are other non constructive results. Peter Mani-Levitska [M] was the first to propose a probabilistic approach: he proved that almost surely any convex body is transformed into a sphere by a Steiner process in which the directions are taken at random with respect to the uniformly distributed probability on the sphere $S^{n-1}\subset \mathbb{R}^n$. This result has been extended to the case of compact seeds [vS] and later to $L_p$ ($1\le p <\infty$) and to continuous functions [V2], where (by a suggestion of A. Burchard) the assumption on the probability on $S^{n-1}$ has been weakened.

Paolo Gronchi constructed a simple but subtle example, which is discussed in [BBGV]. He presented an infinite family of sequences of directions such that the corresponding Steiner processes are not convergent even if the seed $E$ is an ellipse!

This surprising result turned the attention to dense set of directions. In [BKLYZ] it has been showed that if $U\subset S^{n-1}$ is dense and countable and $C$ is a convex body, $U$ can be ordered in such a way that the corresponding Steiner process converges to $C^*$. The result was later improved, showing that there exists a {\it universal reordering} independent from the seed, which can be, more generally, a compact set [V3]. But the proofs are not constructive and there is no idea what this reorderings look like.

To search for other algorithms that make the Steiner process converge, it was natural to look at the well-known class of uniformly distributed sequences, which have been around since more than a century, but were never considered in this context.

The paper [BBGV] contains a first result of this type, proving that the classical Kronecker sequence $\{n \alpha\}$ (where $\alpha$ is an irrational multiple of $\pi$) makes the Steiner process converge for any compact seed $K$.

On the other hand, a subclass of Gronchi's sequences can be used to show that uniform distribution of the sequence is not sufficient for the convergence of Steiner processes.

This observation suggested Conjecture 5.2 which asks if the distinguished class of {\it low-discrepancy} sequences always gives rise to a convergent Steiner process.

The conjecture is still open, but in the meantime a result by R. Asad and A. Burhard (not published but available on Archive [AB]) shows that the well-known low-discrepancy van der Corput binary sequence provides a convergent Steiner process and therefore does not contradict the Conjecture.

In this paper we will prove an analogous result considering the Kakutani-Fibonacci sequence.

\section{Basic definitions}

{\bf 2.1 Steiner symmetrization}\smallskip

Let $M$ be a measurable planar set having finite Lebesgue measure (which will be denoted in this paper by $\lambda$, regardless of the dimension). We will call such sets $L_1$-sets. In the family of $L_1$-sets, we will consider the usual equivalence relation $M_1\sim M_2$ if and only if $\lambda(M_1\triangle M_2)=0$. The quotient space is a metric space with distance $d_1(M_1, M_2)=\lambda(M_1\triangle M_2)$.

If $u\in S^1$ is a unit vector, let us denote by $l_u$ the line through the origin parallel to $u$. 
Denote by $l_{u^{\perp}}$ the line orthogonal to $u$ through the origin. For any $x\in l_{u^{\perp}}$ consider the line $x+l_u$ and intersect it with $M$. Consider now the union of all segments (possibly degenerate and possibly infinite) on $x+l_u $ centered at $x$ whose length equals the exterior measure $\lambda^*(M\cap (x+l_u)$). We will denote this set by $S_uM$ and call it the {\it Steiner symmetral} of $M$ in direction $u$.

Steiner symmetrization has several interesting properties.\smallskip

\noindent 1. $S_uE$ is convex, compact or measurable, if E is convex, compact or measurable, respectively.

\noindent 2. $S_u$ preserves the area and does not increase the perimeter (when it exists). For an update on this topic, see [CCF].

\noindent 3. If $M_1 \subset M_2$, then $S_u M_1 \subset S_u M_2$.

\noindent 4. $d_1 (S_u M_1, S_u M_2)\le d_1 (M_1, M_2)$.

\noindent 5. If $\lim_{n \rightarrow \infty} u_n=u_0$, $\lim_{n \rightarrow \infty} S_{u_n}M=S_{u_0}M$ for any $L_1$-set (Lemma 2.3 of [V2]).
\smallskip

We add to these properties a proposition which is possibly known, but the authors were not able to find a reference.\smallskip

{\bf 2.1.1 Definition}\smallskip

An $L_1$-set $M\subset \mathbb{R}^2$ is called an {\it annulus} centered at $o$ if there exists an $L_1$-set $M_0\subset \mathbb{R}_+$ such that $M=\{(x,y): \sqrt{x^2+y^2}\in M_0\}$, up to a set of measure zero.

Let us denote by $R_{\alpha}E$ the rotation of the set $E$ by the angle $\alpha$ in the anti-clockwise direction.
\bigskip

{\bf 2.1.2 Proposition}\smallskip

Suppose $u,v \in S^1_+$ make an angle $\alpha$ which is an irrational multiple of $\pi$ and $M$ is an $L_1$ set symmetric with respect to $u$ and $v$ up to a set of measure zero. Then $M$ is an annulus centered at the origin.
\smallskip

Proof:

Let $z$ be a point of Lebesgue density $1$ for $M$. This means that given any $\varepsilon >0$ there exists a positive $r_{\varepsilon }$ such that
$$\frac{\lambda(B(z,r_{\varepsilon })\cap M)}{B(z,r_{\varepsilon })}>1-\frac{\varepsilon}{2}\,.$$
Symmetry with respect to $u$ and $v$ implies that $R_{2n\alpha}M=M$ a.e. for each $n\ge 1$.
The points $R_{2n\alpha}z$ are dense on the circle $x^2+y^2=\|z\|^2$.
For each $w$ such that $\|w\|=\|z\|$ there exists an integer $n_{\varepsilon}$ such that $R_{2n_{\varepsilon}\alpha}B(z,r_{\varepsilon })$ and $B(w,r_{\varepsilon })$ are so close to each other that
$$\frac{\lambda(B(w,r_{\varepsilon })\cap M)}{B(w,r_{\varepsilon })}>1-\varepsilon\,.$$
Since $\varepsilon$ is arbitrary, $w$ has also density $1$ with respect to $M$. 

A similar argument shows that if $z$ has density $0$ for $M$, all the points having the same norm also have density $0$. 

This proves that $M_1$, the set of all points which have density $1$ for $M$, is the union of circles centered at $o$. Therefore, $M$ is an annulus.
\bigskip

\noindent {\bf 2.2 Moment of inertia}\smallskip

The moment of inertia of a bounded measurable set $M$ is defined by 
$$\mu(M)=\int_{\mathbb{R}^s}(x^2+y^2))\,dx dxy\,.$$ 
We will use the following properties.\smallskip

\noindent 1. Steiner symmetrization does not increase the moment of inertia, i.e. 
$$\mu(S_uM)\le \mu(M)\,,$$
 with equality if and only if $M$ is symmetric with respect to $u^{\perp}$ up to a set of measure zero.

\noindent 2. If $\cal{F}$ is a bounded family of measurable sets, the functional $(u,M) \longmapsto \mu(S_uM)$ defined on $(S^1_+) \times {\cal F}$ is uniformly continuous. 
\bigskip

\noindent {\bf 2.3 Uniformly distributed sequences of points and discrepancy}\smallskip

A sequence of points $\{x_n\}$ in $[0,1]$  is called {\it uniformly distributed} if for any real-valued continuous function $f$ on $[0,1]$
$$\lim_{N\rightarrow \infty}\frac{1}{N}\sum_{1}^{N}f(x_n)=\int_{0}^{1}f(t)\,dt \,. \eqno(1)$$
Standard references are [KN] and [DT].

The {\it discrepancy} of a sequence $\{x_n\}$ in $[0,1]$ is defined by
$$D_N=\sup _{I\subset [0,1]}\left|\frac {|\{ n \le N, x_n\in I\}|}{N}-\lambda(I)\right|\,,$$ 
the supremum being taken over all the subintervals $I$ of $[0,1]$.
$|\cdot |$ denotes the cardinality.\smallskip

H. Weyl [W] proved that a sequence is uniformly distributed if and only if $\lim_{N\rightarrow \infty} D_N=0$.\smallskip

A sequence $\{x_n\}$ in $[0,1]$ is said to have low-discrepancy if $D_N$ behaves, when $N\rightarrow \infty$, like $\frac{\log N}{N}$. Such sequences provide the best approximation in (1).

The best known examples of of low-discrepancy sequences are the Kronecker sequence $x_ n =\lfloor n \alpha \rfloor$ when $\alpha$ is an irrational number with bounded continued fraction coefficients ([KN], p. 122), the van der Corput sequences of base $b$ ([KN], p. 127), where $b$ is a positive integer, LS-sequences with $L$ and $S$ non negative integers such that $L\ge S$ and $L+S\ge 2$ [C1]. This class contains the van der Corput sequences which correspond to the values $L=b$ and $S=0$. They are low-discrepancy if and only if $L\ge S$. For an update, see [We].\smallskip

The Kakutani-Fibonacci sequence is obtained for $L=S=1$.
\smallskip

The definitions of uniform distribution and discrepancy extend naturally to sequences on $[0,\pi]$ by the map $x\longmapsto \vartheta=\pi x$. It can be extended further to the halfcircle $S^1_+$ (representing the planar directions) by the map $\vartheta \longmapsto e^{i\vartheta}$.
\bigskip

\noindent {\bf 2.4 BV-sets.}\smallskip

In the Section 4 we will work with a class of sets which were introduced by Caccioppoli [Ca1], [Ca2] and used by De Giorgi [DG] in his fundamental paper on the isoperimetric problem: the BV-sets. This is the largest class in which the perimeter can be defined and the isoperimetric problem makes sense.\smallskip

\noindent {\bf Definition}
\smallskip

A bounded set $C$ is a BV-set if its De Giorgi-Caccioppoli perimeter 
$$p(C)= \inf \{\liminf_{k\rightarrow \infty}p(A_k)\}$$
is finite, the greatest lower bound being taken over all sequences $\{A_n\}$ of smooth sets such that $d_1(C,A_n)$ tends to zero and the perimeter of the $A_n$'s is understood in the ordinary sense.
\smallskip

The family $\cal{C}$ of BV-sets satisfies, in particular, the following properties.

1. The Steiner symmetral of a BV-set is a BV-set.

2. The Steiner symmetrization does not increase the perimeter of a BV-set.

3. A family of BV-sets which are contained in a ball $B(o,R)$ and have uniformly bounded perimeters is relatively compact (with respect to $d_1$).

4. The family of BV-sets is dense in $L_1$.\smallskip

An excellent exposition on BV-sets is [T].

\section {\bf  Kakutani-Fibonacci sequence}

In this section we will present the Kakutani-Fibonacci sequence of points and study some of its useful properties. It has been introduced in [C1]. A more explicit algorithm for computing its points can be found in [C2]. This sequence has been studied also from the point of view of dynamical systems [CIV].\smallskip

In [K] Kakutani introduced an elegant and intriguing method to obtain a sequence of {\it partitions} of $[0,1]$.

Fix $0<\alpha <1$ and consider any partition $\delta$ of $[0,1]$. Divide the longest interval(s) of $\delta$ in proportion $\alpha/ (1-\alpha)$ to obtain a finer partition denoted by $\alpha \delta$ (called $\alpha$-refinement of $\delta$). By $\alpha^n \delta$ we denote the $\alpha$-refinement of $\alpha^{n-1}\delta$.
If $\delta=\omega$, the trivial partition of $[0,1]$, the sequence $\{\alpha^n \omega \}$ is called Kakutani's $\alpha$-sequence of partitions.\smallskip
 
 He proved the following theorem.\bigskip
 
 {\bf 3.1 Theorem} \smallskip
 
 The sequence $\{\alpha^n \omega\}$ is uniformly distributed for any $\alpha \in ]0,1[$.
 \smallskip
 
 This means that, given any interval $I\subset [0,1]$ and denoting by $t_n$ the number of intervals of $\alpha^n \omega$ and by $t_n(I)$ the number of those which are contained in $I$, we have
 $$\lim_{n \rightarrow \infty}\frac{t_n(I)}{t_n}=\lambda(I)\,.$$
 
In general it is not simple to track the number of the intervals of $\alpha^n \omega$ and their lengths. There are two interesting exceptions: the cases $\alpha=\frac{1}{2}$ and $\alpha = \frac{1}{2}(\sqrt{5}-1)$, the inverse of the golden ratio, usually denoted in the literature by $\gamma$. In the first case all the intervals of the $n$-th sequence have the same length $\frac{1}{2^n}$. In the second case the intervals of the $n$-th partition have just two lengths: $\gamma^n$ and $\gamma^{n+1}$. This fact follows from the identity $1-\gamma=\gamma^2$ which implies that a long interval is split in any $\gamma$-refinement in the proportion $\gamma/\gamma^2$.

 In [C1] it has been shown that, if we denote by $t_n$, $l_n$ and $s_n$ the number of the intervals, the long intervals and the short intervals of $\gamma^n\omega$, these sequences satisfy the same difference equation $x_{n+2}=x_{n+1}+x_n$, with initial conditions $t_0=1$ and $t_1=2$, $l_0=1$ and $l_1=1$, $s_0=0$ and $s_1=1$. Therefore the three sequences can be expressed with the Fibonacci sequence: $t_n=F_{n+2}$, $l_n=F_{n+1}$ and $s_n=F_n$. This is why this sequence is called the Kakutani-Fibonacci sequence.\smallskip
 
Uniformly distributed sequences of points are far more interesting than uniformly distributed sequences of partitions, so it is worthwhile to note that it is possible to associate efficiently sequences of points to the two sequences of partitions considered above.\smallskip
 
The case $\alpha=\frac{1}{2}$ is well-known ([KN], page 127). Binary fractions can be ordered in the classical van der Corput sequence in base $2$, $\{c_n\}$, which is obtained by means of the {\it radical inverse function} $\Phi$. 

$\Phi$ is defined on the set of all positive integers $\mathbb{N}$. If $\sum_{k=0}^M a_k(n)2^k$ is the binary representation of $n$ (with $M=\lfloor \log_2(n) \rfloor$ and $a_k(n)\in \{0,1\}$), then
 $$c_n=\Phi(n)=\sum_{k=0}^M a_k(n)2^{-k-1}\,.$$
 
 Similarly, the Kakutani-Fibonacci {\it sequence of points} $\{f_n\}$ can be obtained explicitly by ordering the points of the partitions $\gamma^n \omega$.
 
We use a special case of formula (13) of [C2].

This formula is related to R\`eny's representation of real numbers in non-integer base [R]. We will not use results from that deep theory. 

Let us just observe that the equation $1-\gamma=\gamma^2$ shows the ambiguity of the representation, which requires a restriction of the domain of the $\gamma$-radical inverse function we are going to define.
 
Denote by $\mathbb{N}_{\gamma}$ the set of positive integers $n$ which have the property that in their binary representation $a_{k+1}(n)a_{k}(n)=0$ for every $k\ge 0$. 
The {\it $\gamma$-radical inverse function} is defined on $\mathbb{N}_{\gamma}$ as
$$\Phi_{\gamma}(n)=\sum_{k=0}^M a_k(n)\gamma^{k+1}\,. \eqno(2)$$

The first $12$ points of the Kakutani-Fibonacci sequence $\{f_k\}$ are obtained for $n\in \mathbb{N}_{\gamma}\cap \{n: n<2^5\}$ and are, in their order,
$$\gamma,\gamma^2,\gamma^3,\gamma^3+\gamma,\gamma^4,\gamma^4+\gamma, \gamma^4+\gamma^2,\gamma^5,\gamma^5+\gamma,\gamma^5+\gamma^2,\gamma^5+\gamma^3, \gamma^5+\gamma^3+\gamma.$$
To simplify notations, we will denote by $T_k$ the composition $S_{f_k}S_{f_{k-1}}\dots S_{f_1}$.

Let us give a closer look at formula (12) of [12].

It follows from (2) that since $2^{n-1}$ is represented in base $2$ by $100\dots 0$ ($n-1$ zeros), $\Phi_{\gamma}(2^{n-1})=\gamma^n$ is the first element of the sequence $\{f_k\}$ which is an endpoint of $\gamma^n \omega$.
 
We have also $\Phi_{\gamma}(2^{n-1}+1)=\gamma^n+\gamma$, since $2^{n-1}+1$ in base $2$ is represented by $100\dots 1$ (with $n-2$ zeros).

From these observations we can deduce that $\mathbb{N}_{\gamma}$ contains $F_n$ integers $i$ such that $2^{n-1}\le i <2^n-1$.

The two equalities will be used in the next section. Let us rephrase them.

Put $G_k=\sum_{i=1}^k F_i$. Let us recall that $G_k=F_{k+2}-2$. We have that $T_{G_k}=S_{\gamma^k}T_{G_k-1}\,,$ and $T_{G_k+1}=S_{\gamma+\gamma^k}T_{G_k}$.

We will use in the next calculations the $G_k$'s to simplify notation.

\section {\bf  The Kakutani-Fibonacci Steiner process.}

We are now going to prove the main result of the paper, showing that the Steiner process associated to the Kakutani-Fibonacci sequence of directions $\{\vartheta_n=e^{i\pi f_n}\}$ converges with respect to $d_1$ distance if the seed is an $L_1$-set.
\bigskip

Let us begin with a simple lemma.\smallskip

\noindent {\bf 4.1 Lemma.}\smallskip

Let $(X,d)$ be a compact metric space, $\varphi$ a real valued continuous function defined on $X$ and $x_0\in X$ be such that $\varphi(x)>\varphi(x_0)$ for any $x\not=x_0$. If a sequence $\{x_n\}$ is such that $\{\varphi(x_n)\}$ tends to $\varphi(x_0)$, then $\{x_n\}$ tends to $x_0$. \smallskip

Suppose the contrary. Then there exist an $\varepsilon >0$ and a subsequence $\{x_{n_k}\}$ such that $x_{n_k}\not \in B(x_0, \varepsilon)$ for all $k\in \mathbb{N}$. By compactness, it contains a convergent subsequence, which will still be indicated with $\{x_{n_k}\}$. Let $y_0$ be its limit. By construction, $d(y_0,x_0) \ge \varepsilon$ and hence $y_0\not=x_0$. On the other hand $\varphi(y_0)=\lim_{k\rightarrow \infty}\varphi(x_k)=\varphi(x_0)$, a contradiction. \bigskip

\noindent {\bf 4.2 Theorem.}\smallskip

If $M$ is an $L_1$-set and if $\{T_n\}$ is the Kakutani-Fibonacci Steiner process, then
$$\lim_{n\rightarrow \infty}T_nM=M^*$$
 with respect to $d_1$.\smallskip

We first prove the result assuming that $M=C$ is a bounded BV-set.

Consider the subsequence $\{T_{G_k +1}C\}$. We know that if $C\subset B(o,r)$, then $T_{G_k +1}C \subset B(o,r)$. Moreover, the perimeters of the sets of the subsequence are bounded by the perimeter of $C$. Therefore the compactness criterion applies and there exists a convergent subsequence $\{T_{G_{k_j} +1}C\}$. Denote by $K$ its limit.

Since  $T_{G_{k_j} +1}C= S_{G_{k_j} +1}T_{G_{k_j}}C$, the $j$-th term is symmetric with respect to $(e^{i\gamma^{k_j}})^{\perp}$. By the continuity of the Steiner symmetrization, $K$ is symmetric with respect to $(e^0)^{\perp}=(0,1)$.

On the other hand $S_{G_{k_j} +2}C$ is symmetric with respect to $(e^{\gamma^{k_j}+\gamma})^{\perp}$ for each $j$. From the identity 
$$T_{G_{k_j} +2}C= S_{G_{k_j} +2}T_{G_{k_j+1}}C\,,$$
where $S_{G_{k_j} +2}= S_{\gamma^{k_j}+\gamma}$, it follows that $T_{G_{k_j} +2}C$ tends to $S_{\gamma}K$.

The whole sequence $\mu(T_nC)$ tends to $\mu(K)$, therefore
$$\mu(S_{\gamma}K)=\mu(K)\,.$$
This identity implies that $K$ is symmetric with respect to $\gamma^{\perp}$. But $\gamma$ is an irrational multiple of $\pi$ and so is $\gamma^{\perp}$. 

By Proposition 2.1.1 $K$ is, up to a set of measure zero, an annulus. But $S_{\gamma}K=K$ and therefore the intersection of $S_{\gamma}K$ with every line orthogonal to $e^{i\gamma}$ is a segment, hence $K$ is not hollow and it has to be a ball.

In conclusion $K=C^*$ up to a set of measure zero.

As noted, the sequence $\{\mu(T_nC)\}$ is decreasing and converges to $\mu(C^*)$ which is the unique minimizer of $\mu$ on the family $\cal F$ of sets contained in $B(o,r)$ whose measure equals $\lambda(C)$. By Lemma 4.1 the whole sequence $\{T_nC\}$ tends to $C^*$.

Let now be $M$ any $L_1$-set. Fix an arbitrary $\varepsilon >0$. There exists a finite union of rectangles (a bounded BV-set!), denoted by $R$, such that $d_1(M,R)<\frac{\varepsilon}{2}$. 
We may assume that $\lambda(M)=\lambda(R)$.

We have
$$d_1(T_nM, M^*)\le d_1(T_nM,T_nR)+d_1(T_nR, R^*)+d_1(R^*,M^*)\,.
$$
The last term is zero. The first is bounded by $\frac{\varepsilon}{2}$. The second tends to zero and therefore there exists $n_{\varepsilon}$ such that it is bounded by $\frac{\varepsilon}{2}$ for every $n>n_{\varepsilon}$. This proves that $\lim_{n\rightarrow \infty}T_nM=M^*$ in the topology generated by $d_1$.\bigskip

\section {\bf Final remarks.}

In the previous Section we choose to study the behavior of the Kakutani-Fibonacci Steiner process in the class of $L_1$-sets where the topology is defined by the symmetric difference distance $d_1$. 

Other authors treated problems concerning the convergence of Steiner processes in different classes of objects. 

In [AB], the paper whose aim is the most similar to ours, it is proved that the van der Corput Steiner process converges uniformly (to the symmetric decreasing rearrangement) on the class of continuous functions with compact support.

In [BBGV], which motivated this paper, the convergence of the Steiner process is studied in the space of compact sets endowed with the Hausdorff distance. 

In [Kl] Steiner processes are studied in the space of convex bodies with the Hausdorff distance. 

A natural question is, how do these results relate?

The almost complete answer is known [V2]. In Section 4 (Theorem 4.3) of [V2] it is shown that convergence of Steiner processes in the space of $L_1$-sets implies the convergence of Steiner processes in the space of $L_p$ functions (for $p\ge 1$) to the symmetric decreasing rearrangement in the $\|.\|_p$ norm. An almost immediate consequence is Theorem 4.4, which states that if the Steiner process converges in $\|.
\cdot \|_1$ norm to the symmetric decreasing rearrangement for every function in $L_1$, it converges uniformly to the symmetric decreasing rearrangement if the seed is a continuous function with compact support. 

Finally Theorem 5.1 in Section 5 asserts that uniform convergence of the Steiner process for continuous function with compact support to the symmetric decreasing rearrangement implies its convergence in the Hausdorff distance to $K^*$ for every compact seed $K$.

To close the circle we only need to prove that convergence of a Steiner process applied to a compact set $K$ to $K^*$ in the Hausdorff distance implies its convergence to $M^*$ for any $L_1$-set $M$.\smallskip

We need the following lemma. It is known, but we will prove it to show how powerful tool is the Steiner symmetrization.\smallskip

\noindent {\bf Lemma 5.1}\smallskip

If $A, B \subset \mathbb{R}^m$ are two $L_1$-sets, and $A^*$ and $B^*$ are the balls centered at $o$ having the same measure as $A$ and $B$ respectively, then
$$d_1(A^*,B^*)\le d_1(A,B)\,.$$

Let $\{T_n=S_nS_{n-1}\dots S_1\}$ be any Steiner process which converges to $M^*$ for every $L_1$-set $M$. By one of the properties of the Steiner symmetrization
$$d_1(T_nA,T_nB)\le d_1(A,B)$$
for every $n$. Letting $n\rightarrow \infty$ we get the conclusion.
\smallskip

\noindent {\bf Lemma 5.2}\smallskip

If $K\subset \mathbb{R}^m$ is a compact set and $\{T_n\}$ is a Steiner process such that the sequence $\{T_nK\}$ converges in the Hausdorff distance to a compact set $L$, then 
$$\lim_{n\rightarrow \infty}\lambda(T_nK\triangle L)=0\,.$$
This is Lemma 3.1, claim (ii), of [BBGV]. \smallskip

\noindent {\bf Proposition 5.3}\smallskip

Suppose that a Steiner symmetrization converges (in the Hausdorff distance) for any compact seed. Then it converges (with respect to $d_1$) for any $L_1$-seed $M$.\smallskip

Let $M$ be a $L_1$-set. Then, for every $\varepsilon>0$ there exists a compact set $K_{\varepsilon}\subset M$ such that $\lambda(M\setminus K_{\varepsilon})<\frac{\varepsilon}{3}$. We know that $T_nK_{\varepsilon} \subset T_n M$ for every $n\ge 1$. We have
$$d_1(T_nM,M^*)\le d_1(T_nM,T_nK_{\varepsilon})+d_1(T_nK_{\varepsilon}, K_{\varepsilon}^*)+d_1(K_{\varepsilon}^*,M^*)\,.$$
The first term is bounded by $d_1(M, K_{\varepsilon})$ and hence by $\frac{\varepsilon}{3}$. The same bound holds for the third term by Lemma 5.1. The middle term tends to zero by Lemma 5.2. Therefore there exists $n_{\varepsilon}$ such that for $n\ge n_{\varepsilon}$ it is smaller then $\frac{\varepsilon}{3}$, and the desired conclusion follows.
\smallskip

Steiner processes have been studied first in the family of convex bodies. Since a convex body is  a compact set, any result proved in one of the settings mentioned above holds also for convex bodies. 

The Open Problem 6.3 of [BBGV] asks the following question: {\it Assume that a sequence of directions $\{u_n\}$ is such that $\{S_{u_n}\dots S_1C\}$ converges to $C^*$ for each convex body $C$. Is it true that $\{S_{u_n}\dots S_1K\}$ converges to $K^*$ for each compact set $K$?}

The authors do not know of any result which holds for convergence of Steiner processes on convex bodies which does not hold for compact seeds. But there does not exist (so far, at least) a general and short way to extend a result which holds for convex bodies to the class of compact sets.

As an example of the difficulties addressed by Problem 6.3, let us  mention Mani-Levitska's question [M], posed in 1986. He proved his result when the seed is a convex body and asked if it holds also for compact seeds. The positive answer has been proved only in 2006 [vS]. Another solution has been found in 2013 [V2]. The two papers arrived at the desired result using completely different methods, both after a long detour. 
\bigskip

Problem 5.2 remains unanswered. This paper partially confirms it, adding another example of low-discrepancy sequence of directions which makes the Steiner process convergent.
\bigskip\bigskip

 \noindent {\bf References}
 \bigskip
 
\noindent [AB]  R. Asad, A. Burchard, A.: {\it Steiner symmetrization along a certain equidistributed sequence of directions}, posted on arXiv:2005.13597 May 2020.\bigskip

\noindent [BBGV] G. Bianchi, A. Burchard, P. Gronchi, A. Vol\v{c}i\v{c}, {\it Convergence in shape of Steiner Symmetrizations}, Indiana Univ. Math. J. {\bf 61}, no. 4 (2012), 1695-1710 \bigskip

\noindent [Ca1] R. Caccioppoli, {\it Misura e integrazione sugli insiemi dimensionalmente orientati I}, Rend. Accad. Naz. Lincei Cl. Sci. Fis. Mat. Nat., {\bf 12} (1952), 3-11  \bigskip

\noindent [Ca2] R. Caccioppoli, {\it Misura e integrazione sugli insiemi dimensionalmente orientati II}, Rend. Accad. Naz. Lincei Cl. Sci. Fis. Mat. Nat., {\bf 12} (1952), 137-146 \bigskip

\noindent [C1] I. Carbone, {\it Discrepancy of LS-sequences of partitions and points}, Annali di Matematica Pura e Applicata, ({\bf 4}), {\bf 191(4)} (2012), 819-844 \bigskip

\noindent [C2] I. Carbone, {\it Extension of van der Corput algorithm to LS-sequences}, Applied Mathematics and Computation, {\bf 255} (2015), 207-2013
\bigskip

\noindent [CCF] A. Cianchi, M. Chleb\'{\i}k, N. Fusco, {\it The perimeter inequality under Steiner
symmetrization: Cases of equality}, Annals of Mathematics, {\bf 162} (2005), 525-555 
\bigskip

\noindent [CIV] I. Carbone, M.R. Iac\`o, A. Vol\v{c}i\v{c}, {\it A dynamical system approach to the Kakutani-Fibonacci sequence}, Ergodic Theory and Dynamical Systems, {\bf 34}(6) (2014), 1794-1806, 

\bigskip

\noindent [DG] E. De Giorgi, {\it Su una teoria generale della misura $(r-1)$-dimensionale in uno spazio a $r$ dimensioni}, Annali Mat. Pura e Appl. (4) {\bf 36} (1954), 191-213
\bigskip

\noindent  [DT] M. Drmota,  R.F. Tichy, {\it  Sequences, discrepancies and applications}, (1997) Springer LNM {\bf 1651} 
\bigskip

\noindent [E] F. Edler: {\it Vervollst\"andigung der Steinerschen elementargeometrichen Beweise f\"ur den Satz, das der Kreis gr\"osseren Fl\"acheninhalt besitz als jede andere ebene Figur gleich grossen Umfangs}, Nachr. Ges. Wiss. G\"ottingen, (1882), 73-80
 \bigskip
 
\noindent [K] S. Kakutani, S. {\it A problem of equidistribution on the unit
interval $[0,1]$.} Measure theory (Proc. Conf., Oberwolfach, 1975),  369-375. Springer LNM, {\bf 541} (1976)
\bigskip

\noindent [Kl] D. Klein, {\it Steiner symmetrization using a finite set of directions}, Adv. in Appl. Math. {\bf 48} (2012), 1322-1338 
\bigskip

\noindent [KN] L. Kuipers, H. Niederreiter, {\it Uniform distribution of sequences}, Wiley (1974)
\bigskip

\noindent [M] P. Mani-Levitska, {\it Random Steiner symmetrizations}, Stud. Sci. Math. Hung. {\bf 21}(3-4),  (1986), 373-378
\bigskip

\noindent [R] A. R\'enyi, {\it Representations for real numbers and their ergodic properties}, Acta Mathematica Academiae Scientiarum Hungaricae, {\bf 8} (3-4) (1957), 477-493
\bigskip

\noindent [T] G. Talenti, {\it The standard isoperimetric theorem}. In: Gruber, P.M., Wills, J.M. (eds.) Handbook of convex geometry, 73-123, North-Holland, (1993)
\bigskip

\noindent [V1] A. Vol\v{c}i\v{c}, {\it A generalization of Kakutani`s splitting procedure}, Ann. Mat.
Pura e Appl. (4) {\bf 190} no. 1, (2011), 45-54  \bigskip

 \noindent [V2] A. Vol\v{c}i\v{c}, {\it Random Steiner symmetrizations of sets and functions}, Calc.
Var. Partial Differential Equations {\bf 46}, no. 3-4, (2013), 555-569, \bigskip

\noindent [V3] A. Vol\v{c}i\v{c}, {\it Iterations of random symmetrizations}, Annali di Mat. Pura e
Appl., {\bf 195} n. 5, (2016), 1685-1692  \bigskip

\noindent [We] C. Wei\ss, {\it On the classification of LS-sequences}, Unif. Distrib. Theory {\bf 13} , no. 2, (2018), 83-92
\bigskip

\noindent [W] H. Weyl, {\it \"Uber ein Problem aus dem Gebiete der diophantischen Approximationen}, Nach. Ges. Wiss. G\"ottingen, Math.-phys. Kl.  (1914), 234-244  \bigskip

\end{document}